\documentclass[12pt]{article}

\usepackage{amsmath}
\usepackage{amssymb}
\usepackage[dvips]{graphicx}

\begin{document}
\title{{\bf 
{\large Prime component-preservingly amphicheiral link
with odd minimal crossing number}}
\footnotetext[0]{%
2010 {\it Mathematics Subject Classification}: 
57M25, 57M27.\\
{\it Keywords}:
component-preservingly amphicheiral link; minimal crossing number; 
Tait's conjecture ; invertibility.
}}
\author{
{\footnotesize 
Teruhisa KADOKAMI and Yoji KOBATAKE}}
\date{{\footnotesize March 11, 2015}}
\maketitle

\newcommand{\circlenum}[1]{{\ooalign{%
\hfill$\scriptstyle#1$\hfill\crcr$\bigcirc$}}}

\newcommand{\svline}[1]{\multicolumn{1}{|c}{#1}}
\newfont{\bg}{cmr10 scaled\magstep4}
\newcommand{\bigzerol}{\smash{\hbox{\bg 0}}}
\newcommand{\bigzerou}{\smash{\lower1.7ex\hbox{\bg 0}}}

\newcommand{\bsquare}{\hbox{\rule{6pt}{6pt}}}
\newcommand{\qed}{\hbox{\rule[-2pt]{3pt}{6pt}}}
\newcommand{\Int}{\mathrm{Int}\ \! }
\newcommand{\Ker}{\mathrm{Ker}\ \! }
\newcommand{\Ig}{\mathrm{Im}\ \! }
\newcommand{\aug}{\mathrm{aug}\ \!}
\newcommand{\pj}{\mathrm{pr}\ \!}
\newcommand{\Tor}{\mathrm{Tor}\ \!}
\newcommand{\Spin}{\mathrm{Spin}\ \!}
\newcommand{\Eul}{\mathrm{Eul}\ \!}
\newcommand{\Vect}{\mathrm{Vect}\ \!}
\newcommand{\HULL}{\mathrm{HULL}\ \!}
\newcommand{\real}{\mathrm{real}\ \!}
\newcommand{\rank}{\mathrm{rank}\ \!}
\newcommand{\ord}{\mathrm{ord}\ \!}
\newcommand{\Sign}{\mathrm{Sign}\ \!}
\newcommand{\Hom}{\mathrm{Hom}\ \!}
\newcommand{\ad}{\mathrm{ad}\ \!}
\newcommand{\Det}{\mathrm{Det}\ \!}
\newcommand{\lk}{\mathrm{lk}\ \!}
\newcommand{\pt}{\mathrm{pt}}
\newcommand{\al}{$\alpha$}
\newcommand{\dis}{\displaystyle}

\newtheorem{df}{Definition}[section]
\newtheorem{lm}[df]{Lemma}
\newtheorem{theo}[df]{Theorem}
\newtheorem{re}[df]{Remark}
\newtheorem{pr}[df]{Proposition}
\newtheorem{ex}[df]{Example}
\newtheorem{co}[df]{Corollary}
\newtheorem{cl}[df]{Claim}
\newtheorem{qu}[df]{Question}
\newtheorem{pb}[df]{Problem}
\newtheorem{cj}[df]{Conjecture}

\makeatletter
\renewcommand{\theequation}{%
\thesection.\arabic{equation}}
\@addtoreset{equation}{section}
\makeatother

\begin{abstract}
{\footnotesize 
\setlength{\baselineskip}{10pt}
\setlength{\oddsidemargin}{0.25in}
\setlength{\evensidemargin}{0.25in}
For every odd integer $c\ge 21$,
we raise an example of 
a prime component-preservingly amphicheiral link
with the minimal crossing number $c$.
The link has two components, 
and consists of an unknot and a knot 
which is $(-)$-amphicheiral with odd minimal crossing number.
We call the latter knot a {\it Stoimenow knot}.
We also show that the Stoimenow knot is not invertible
by the Alexander polynomials.}
\end{abstract}

\tableofcontents

\section{Introduction}\label{sec:intro}
Let $L=K_1\cup \cdots \cup K_r$ be 
an oriented $r$-component link in $S^3$.
A $1$-component link is called a knot.
For an oriented knot $K$, we denote the orientation-reversed knot by $-K$.
If $\varphi$ is an orientation-reversing homeomorphism of $S^3$ so that
$\varphi(K_i)=\varepsilon_{\sigma(i)} K_{\sigma(i)}$ for all $i=1, \ldots, r$
where $\varepsilon_i=+$ or $-$,
and $\sigma$ is a permutation of $\{1, 2, \ldots, r\}$, 
then $L$ is called
an {\it $(\varepsilon_1, \ldots, \varepsilon_r; \sigma)$-amphicheiral link}.
A term ``amphicheiral link"  is used as a general term for
an $(\varepsilon_1, \ldots, \varepsilon_r; \sigma)$-amphicheiral link.
If $\varphi$ can be taken as an involution (i.e.\ $\varphi^2=\mathrm{id}$),
then $L$ is called a {\it strongly} amphicheiral link.
If $\sigma$ is the identity, then an amphicheiral link is called
a {\it component-preservingly amphicheiral link}, and
$\sigma$ may be omitted from the notation.
If every $\varepsilon_i=\varepsilon$ is identical for all $i=1, \ldots, r$
(including the case that $\sigma$ is not the identity), then 
an $(\varepsilon_1, \ldots, \varepsilon_r; \sigma)$-amphicheiral link
is called an $(\varepsilon)$-amphicheiral link.
We use the notations $+=+1=1$ and $-=-1$.
For the case of invertibility, we only replace $\varphi$ 
with an orientation-preserving homeomorphism of $S^3$.
We refer the reader to \cite{Wh, Hi, Kd1, Kd2, Kd3, KK}.

\medskip

The minimal crossing number of an alternating amphicheiral link
is known to be even (cf.\ \cite[Lemma 1.4]{Kd3})
from the positive answer for the {\it flyping conjecture}
due to W.~Menasco and M.~Thistlethwaite \cite{MT}.
The flyping conjecture is one of famous Tait's conjectures
on alternating links,
and it is also called Tait's conjecture III in \cite{St1}.
The positive answer for the flyping conjecture implies 
those of Tait's conjecture I on the minimal crossing number
(cf.\ \cite{Mu1}), 
and Tait's conjecture II on the writhe (cf.\ \cite{Mu2}).
A.~Stoimenow \cite[Conjecture 2.4]{St1} sets a conjecture:
\begin{center}
``{\it Amphicheiral (alternating?) knots have even crossing number}."
\end{center}
as Tait's conjecture IV
by guessing what Tait had in mind (i.e.\ Tait has not stated it explicitly).
We pose the following conjecture:
\begin{cj}\label{cj:even}
{\rm (a generalized version of Tait's conjecture IV)}
The minimal crossing number of an amphicheiral link is even.
\end{cj}
For the case of alternating amphicheiral links,
Conjecture \ref{cj:even} is affirmative
as mentioned above from the answer for Tait's conjecture II.
Hence it motivates to find an amphicheiral link with
odd minimal crossing number.
If there exists a counter-example for Conjecture \ref{cj:even}, 
then it should be non-alternating.

\medskip

A non-split link is {\it prime} if it is not a connected sum of non-trivial links.
We assume that a prime link is non-split.
There exists a prime amphicheiral knot 
with minimal crossing number $15$ in the table of 
J.~Hoste, M.~Thistlethwaite and J.~Weeks \cite{HTW},
which gives a negative answer for Conjecture \ref{cj:even}
(the original Tait's conjecture IV).
The knot is named $15_{224980}$ (Figure 1).
\begin{figure}[htbp]
\begin{center}
\includegraphics[scale=0.5]{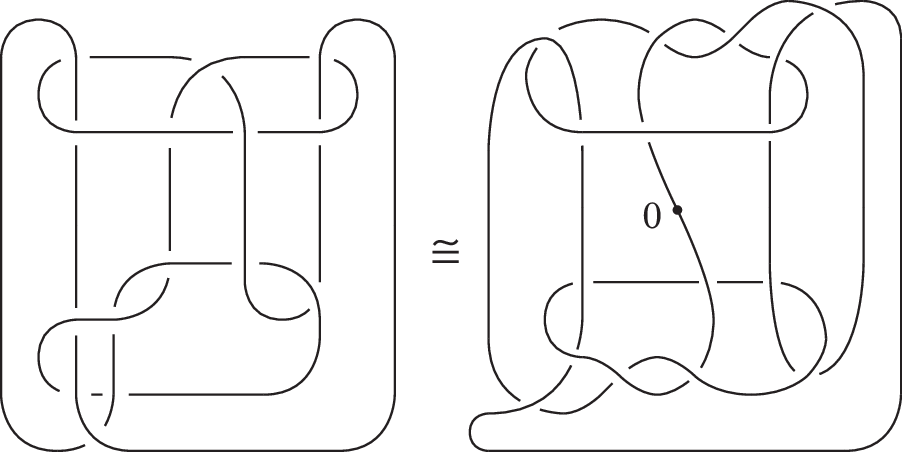} 
\label{fig:HTW}
\caption{$15_{224980}$}
\end{center}
\end{figure}
Stoimenow \cite{St2} showed that
for every odd integer $c\ge 15$,
there exists an example of a prime amphicheiral knot
with minimal crossing number $c$.
The case $c=15$ corresponds to $15_{224980}$.
We call the sequence of knots {\it Stoimenow knots} (see Section \ref{sec:Stoimenow}).
He also pointed out that
there are no such examples for the case $c\le 13$.

\medskip

The first author and A.~Kawauchi \cite{KK},
and the first author \cite{Kd3} determined 
prime amphicheiral links with 
minimal crossing number up to $11$.
Then there are two prime amphicheiral links
with odd minimal crossing numbers
named $9_{61}^2$ and $11_{n247}^2$ (Figure 2),
where we use modified notations
from Rolfsen's table \cite{Ro} and
Thistlethwaite's table
on the web site maintained by
D.~Bar-Natan and S.~Morrison \cite{BM}.
\begin{figure}[htbp]
\begin{center}
\includegraphics[scale=0.6]{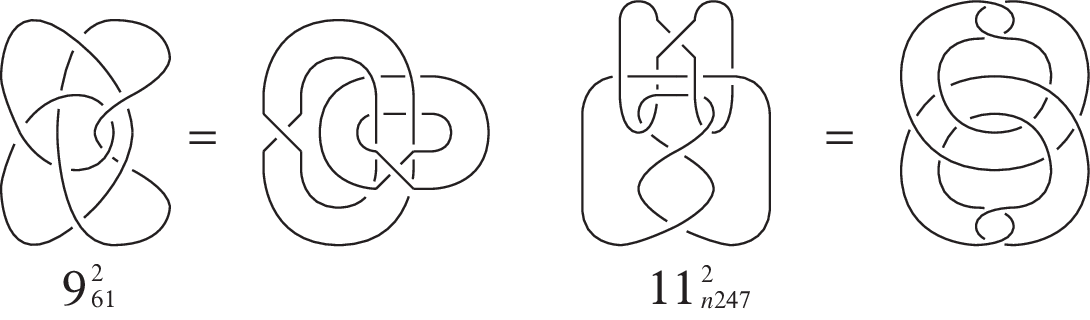} 
\label{fig:odd}
\caption{$9_{61}^2$ and $11_{n247}^2$}
\end{center}
\end{figure}
These examples show that
Conjecture \ref{cj:even} is negative for links.
Since both $9_{61}^2$ and $11_{n247}^2$
are not component-preservingly amphicheiral,
we ask the following question (see also Question \ref{qu:mincr}) :
\begin{qu}\label{qu:cp}
Is there a prime component-preservingly amphicheiral link 
with odd minimal crossing number ?
\end{qu}
If we remove `prime' from Question \ref{qu:cp},
then we can obtain nugatory examples by 
taking split sum of a Stoimenow knot and an unknot,
or connected sum of Stoimenow knot and the Hopf link.
Our main theorem is an affirmative answer for
Question \ref{qu:cp} which is a negative answer
for Conjecture \ref{cj:even}:
\begin{theo}\label{thm:main}
For every odd integer $c\ge 21$,
there exists
a prime component-preservingly amphicheiral link
with minimal crossing number $c$ (Figure 10).
\end{theo}
Our example is a $2$-component link with linking number $3$
whose components are a Stoimenow knot and an unknot.
We prove it in Section \ref{sec:proof}.
The proof is divided into three parts such as
to show amphicheirality,
to determine the minimal crossing number, and
to show primeness.
We can immediately see its amphicheirality by construction.
Though to find the way of linking of the two components
was not so easy,
to determine the minimal crossing number is easy
by the help of Stoimenow's result \cite{St2}
(cf.\ Theorem \ref{thm:Stoimenow}).
In \cite{St2}, to determine the minimal crossing number
and to show primeness of his knot were very hard.
Finally we show primeness by using the Kauffman bracket
(cf.\ Subsection \ref{ssec:Kauffman}).
This part is also eased by Stoimenow's result.
In Section  \ref{sec:noninv}, 
by R.~Hartley \cite{Ha}, R.~Hartley and A.~Kawauchi \cite{HK}, 
and A.~Kawauchi \cite{Kw1}'s
necessary conditions on the Alexander polynomials of amphicheiral knots,
we show that a Stoimenow knot is not invertible (Theorem \ref{th:Stoimenow}).

\section{Link invariants}\label{sec:inv}

\subsection{Kauffman bracket}\label{ssec:Kauffman}

Let $L$ be an $r$-component oriented link, 
and $D$ a diagram of $L$.
Firstly we regard $D$ as an unoriented diagram.
On a crossing of $D$, a {\it splice} is a replacement from 
the left-hand side (the crossing) to the right-hand side as in Figure 3.
Precisely, a {\it $0$-splice} is to the upper right-hand side, and
an {\it $\infty$-splice} is to the down right-hand side, respectively.
\begin{figure}[htbp]
\begin{center}
\includegraphics[scale=1.5]{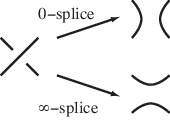} 
\label{splice}
\caption{splice}
\end{center}
\end{figure}
The resulting diagram is a {\it state}, and
it is a diagram of an unlink without crossings.
Let $s$ be a state, 
$|s|$ the number of components of $s$,
$t_0(s)$ the number of $0$-splices to obtain $s$,
$t_{\infty}(s)$ the number of $\infty$-splices to obtain $s$,
$t(s)=t_0(s)-t_{\infty}(s)$,
and $\mathcal{S}$ the set of states from $D$.
Let $A$ be an indeterminate, and $d=-A^2-A^{-2}$.
Then
$$\langle D\rangle =\sum_{s\in \mathcal{S}}
A^{t(s)}d^{|s|-1}\in \mathbb{Z}[A, A^{-1}]$$
is the {\it Kauffman bracket} of $D$, and
\begin{equation}\label{eq:f}
f_L(A)=(-A^3)^{-w(D)}\langle D\rangle
\end{equation}
is the {\it $f$-polynomial} of $L$ 
where $w(D)$ is the writhe of $D$ as an oriented diagram.
Then $f_L(A)$ is an invariant of $L$, and
\begin{equation}\label{eq:J}
V_L(t)=f_L\left(t^{\frac 14}\right)\in 
\mathbb{Z}\left[t^{\frac 12}, t^{-\frac 12}\right]
\end{equation}
is the {\it Jones polynomial} of $L$.
We denote $\langle D\rangle$ as $\langle D\rangle(A)$
when we emphasis it as a function of $A$.
We have the following facts:
\begin{lm}\label{lm:Kauffman}
Let $L$ be an $r$-component oriented link, and $D$ a diagram of $L$.
\begin{enumerate}
\item[(1)]
The Kauffman bracket $\langle D\rangle$ is 
an invariant of $L$ up to multiplications of $(-A^3)$.
In particular, if we substitute a root of unity for $A$
and take its absolute value, then
it is an invariant of $L$, which is a non-negative real number.

\item[(2)]
We have the following skein relation (Figure 4)
which can be an axiom of the Kauffman bracket:
\begin{figure}[htbp]
\begin{center}
\includegraphics[scale=0.8]{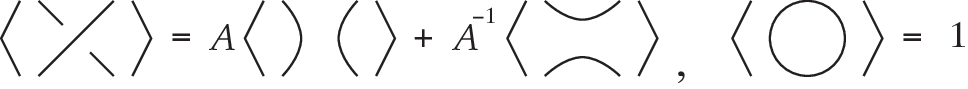} 
\label{skein1}
\caption{skein relation I}
\end{center}
\end{figure}

\item[(3)]
Let $L_i\ (i=1, 2)$ be a link,
$D_i$ a link diagram of $L_i$,
and $D_1\amalg D_2$ ($L_1\amalg L_2$, respectively)
the split sum of $D_1$ and $D_2$ ($L_1$ and $L_2$, respectively).
Then we have
$$\langle D_1\amalg D_2\rangle
=d\langle D_1\rangle\langle D_2\rangle,\ 
f_{L_1\amalg L_2}(A)=d\cdot f_{L_1}(A)f_{L_2}(A).$$

\item[(4)]
Let $L_i\ (i=1, 2)$ be a link,
$D_i$ a link diagram of $L_i$,
and $D_1\sharp D_2$ ($L_1\sharp L_2$, respectively)
the connected sum of $D_1$ and $D_2$ ($L_1$ and $L_2$, respectively).
Then we have
$$\langle D_1\sharp D_2\rangle
=\langle D_1\rangle\langle D_2\rangle,\ 
f_{L_1\sharp L_2}(A)=f_{L_1}(A)f_{L_2}(A).$$

\item[(5)]
We have a skein relation as in Figure 5:
\begin{figure}[htbp]
\begin{center}
\includegraphics[scale=0.8]{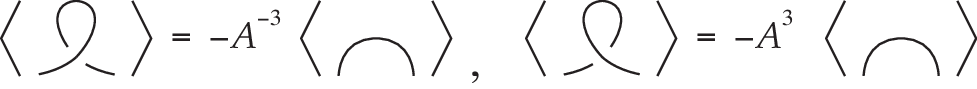} 
\label{skein2}
\caption{skein relation II}
\end{center}
\end{figure}

\item[(6)]
Let $D^*$ ($L^*$, respectively) be the mirror image of $D$ ($L$, respectively).
Then we have
$$\langle D^*\rangle(A)
=\langle D\rangle(A^{-1}),\ 
f_{L^*}(A)=f_L(A^{-1}).$$

\item[(7)]
$f_L(A)\in A^{2(r+1)}\cdot \mathbb{Z}[A^4, A^{-4}]$.

\item[(8)]
Let $\zeta$ be a primitive $8$-th root of unity
(i.e.\ $\zeta^4=-1$ and $\zeta^8=1$).
Suppose that the number of the crossing number of $D$ is even .
Then $\langle D\rangle (\zeta)$
is an integer or of the form $\sqrt{-1}\times (\mbox{integer})$,
which depends on $r$ and the writhe.
In particular, for $r=1$, 
$\langle D\rangle (\zeta)$ is an integer
if and only if the writhe is $0\ (\mathrm{mod}\ \! 4)$.

\item[(9)]
Let $\zeta$ be a primitive $8$-th root of unity.
Then we have $|\langle D\rangle (\zeta)|=|V_L(-1)|$.

\end{enumerate}
\end{lm}

Lemma \ref{lm:Kauffman} (8) is obtained from (7) and (\ref{eq:f}), and
it is a special case of (1).
Lemma \ref{lm:Kauffman} (9) is obtained from (\ref{eq:J}).

\medskip

Let $T_m$ be an $m$-half twist tangle for $m\in \mathbb{Z}$,
and $T_{\infty}$ a tangle in Figure 6.
\begin{figure}[htbp]
\begin{center}
\includegraphics[scale=0.8]{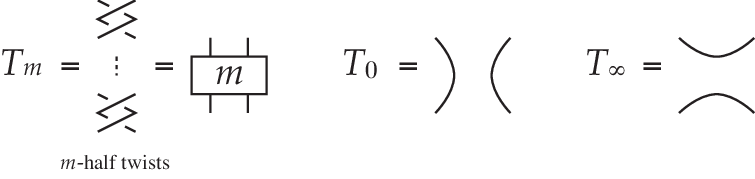} 
\label{mfwist}
\caption{$m$-half twists}
\end{center}
\end{figure}

By Lemma \ref{lm:Kauffman} (2), (3), (4) and (5),
we have the following:
\begin{lm}\label{lm:Tm}
\begin{enumerate}
\item[(1)]
We have
$$\langle T_m\rangle=A^m\langle T_0\rangle
+\alpha_m(A)
\langle T_{\infty}\rangle,$$
where
$$\alpha_m(A)
=A^{m-2}\cdot \frac{1-(-A^{-4})^m}{1-(-A^{-4})}.$$

\item[(2)]
$\alpha_{-m}(A)=\alpha_m(A^{-1}).$

\item[(3)]
Let $\zeta$ be a primitive $8$-th root of unity.
Then we have
$$\alpha_m(\zeta)=m\zeta^{m-2}\quad \mbox{and}\quad
\alpha_m(\zeta)\cdot \alpha_{-m}(\zeta)=m^2.$$

\end{enumerate}
\end{lm}

\subsection{Alexander and Conway polynomials}\label{ssec:Alexander}

Let $L$ be an oriented link, and $D$ a diagram of $L$.
Pick a crossing $c$ of $D$.
If $c$ is a positive crossing (a negative crossing, respectively), 
then we denote $D$ by $L_+$ ($L_-$, respectively).
If $c$ is smoothed with preserving the orientation, 
then we denote $D$ by $L_0$.
We call a pair $(L_+, L_-, L_0)$ a {\it skein triple} (Figure 7).

\begin{figure}[htbp]
\begin{center}
\includegraphics[scale=0.8]{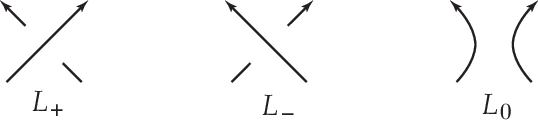} 
\label{skein3}
\caption{skein triple}
\end{center}
\end{figure}

For an oriented link $L$, the {\it Conway polynomial} of $L$ is denoted by $\nabla_L(z)$
which is an element of $\mathbb{Z}[z]$.
For a skein triple $(L_+, L_-, L_0)$, 
the Conway polynomial is defined by the following skein relation : 
$$\nabla_{L_+}(z)-\nabla_{L_-}(z)=z\nabla_{L_0}(z),\quad 
\nabla_O(z)=1,$$
where $O$ is the trivial knot.

\begin{lm}\label{lm:Conway}
Let $L$ be an $r$-component oriented link, and $L^*$ the mirror image of $L$.
Then we have
$$\nabla_{L^*}(z)=\nabla_L(-z).$$
More precisely, 
$\nabla_{L^*}(z)=\nabla_L(z)$ if $r$ is odd, and 
$\nabla_{L^*}(z)=-\nabla_L(z)$ if $r$ is even.
\end{lm}

For an $r$-component oriented link $L$, 
the {\it (normalized one variable) Alexander polynomial} 
${\mit \Delta}_L(t)$ is defined by
$${\mit \Delta}_L(t)=\nabla_L\left(t^{\frac 12}-t^{-\frac 12}\right)
\in \mathbb{Z}\left[ t^{\frac 12}, t^{-\frac 12}\right].$$

For $A, B\in \mathbb{Z}\left[ t^{\frac 12}, t^{-\frac 12}\right]$, 
$A\doteq B$ implies 
$A=\pm t^{\frac m2}B$ for some $m\in \mathbb{Z}$.
For $f, g\in \mathbb{Z}[z]$ or $\mathbb{Z}\left[ t^{\frac 12}, t^{-\frac 12}\right]$, 
if they are equal as elements in 
$(\mathbb{Z}/d\mathbb{Z})[z]$ or $(\mathbb{Z}/d\mathbb{Z})\left[ t^{\frac 12}, t^{-\frac 12}\right]$, 
then we denote by $f=_dg$.
For an oriented link $L$, if $\nabla_L(z)$ and ${\mit \Delta}_L(t)$ are regarded as 
elements in $(\mathbb{Z}/d\mathbb{Z})[z]$ and 
$(\mathbb{Z}/d\mathbb{Z})\left[ t^{\frac 12}, t^{-\frac 12}\right]$ respectively, 
then we call them the {\it mod $d$ Conway polynomial} of $L$ and 
the {\it mod $d$ Alexander polynomial} of $L$ respectively.

\section{Stoimenow knots}\label{sec:Stoimenow}

Let $\sigma_i$\ $(i=1, \ldots, m-1)$ be
a generator of the $m$-string braid group,
and $\delta_i$ and $\overline{\delta}_i$\ $(i=1, \ldots, m-1)$
tangles in Figure 8.
\begin{figure}[htbp]
\begin{center}
\includegraphics[scale=0.8]{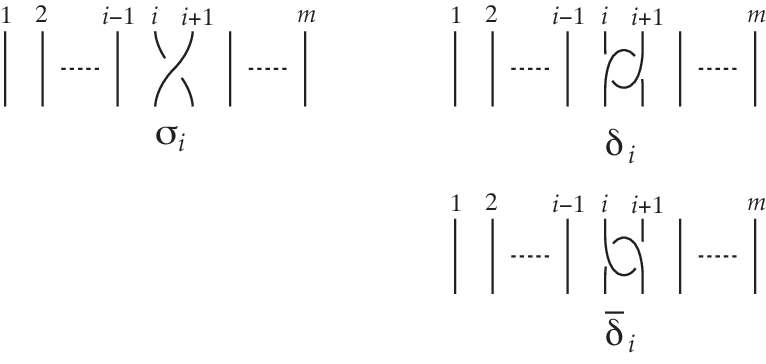} 
\label{braid}
\caption{generator $\sigma_i$ of braid group, 
and $\delta_i$ and $\overline{\delta}_i$}
\end{center}
\end{figure}
For an odd number $n\ge 15$,
a {\it Stoimenow knot} with crossing number $n$,
denoted by $S_n$, is the closure of the following composition
of $\sigma_i$, $\delta_i$ and $\overline{\delta}_i$\ $(i=1, \ldots, m-1)$:
$$\begin{array}{ll}
3\ -1\ 2^2\ 3^{2k}\ 4\ -3\ 2\ -1\ (-2)^{2k}\ (-3)^2\ 4\ -2
& (n=4k+11),\medskip \\
\delta_3\ -1\ 2^2\ 3^{2k}\ 4\ -3\ 2\ -1\ (-2)^{2k}\ (-3)^2\ 4\ 
\overline{\delta}_2
& (n=4k+13),
\end{array}$$
where in the sequence above, $m=5$, 
$\sigma_i$ is translated into $i$ and
$\sigma_i^{-1}$ is translated into $-i$, and
$i^l$ implies that $i$ is repeated $l$ times with $l\ge 1$.
The former is {\it of type I}, and
the latter is {\it of type II}, respectively.
Note that $S_{15}=15_{224980}$ in Figure 1,
and both two tangles above have $(n+1)$ crossings.
We can see strong $(-)$-amphicheirality of $S_n$ from
its diagram with $(n+1)$ crossings in the righthand side of Figure 9.
\begin{theo}\label{thm:Stoimenow}
{\rm (Stoimenow \cite{St1, St2})}
A Stoimenow knot $S_n$ is a prime strongly $(-)$-amphicheiral knot
with minimal crossing number $n$.
\end{theo}

\begin{figure}[htbp]
\begin{center}
\includegraphics[scale=0.6]{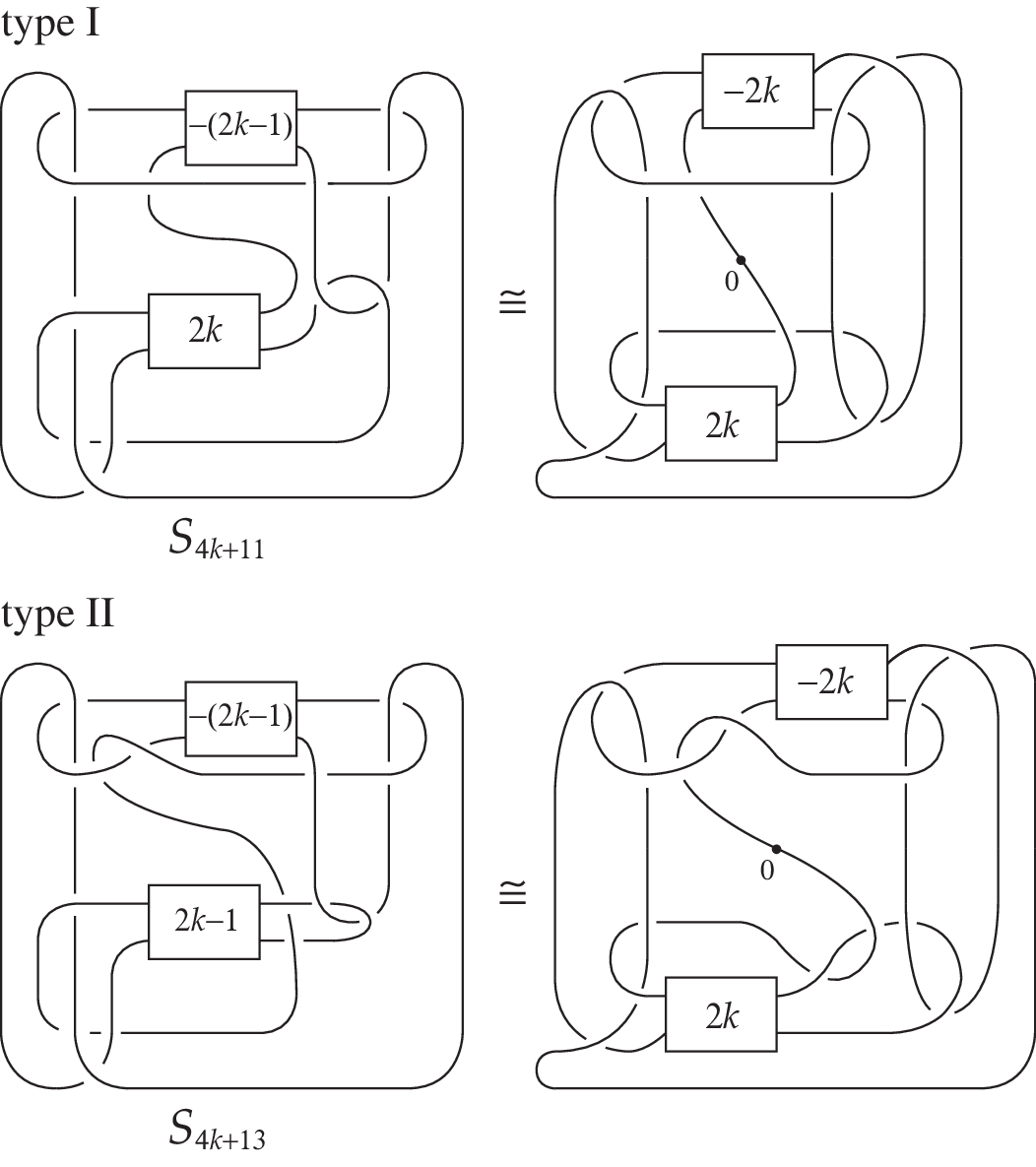} 
\label{Stoimenow}
\caption{Stoimenow knot $S_n$}
\end{center}
\end{figure}

\newpage

\section{Proof of Theorem \ref{thm:main}}\label{sec:proof}

We take a $2$-component link $L_n=S_n\cup U$ whose components are
a Stoimenow knot $S_n$ and an unknot $U$ as in Figure 10.
The link $L_n$ is {\it of type I} if $S_n$ is of type I, and
is {\it of type II} if $S_n$ is of type II.
We prove that
$L_n$ is a prime component-preservingly amphicheiral link
with minimal crossing number $n+6$,
where $n+6$ is odd with $n+6\ge 21$
because $n$ is odd with $n\ge 15$.

\begin{figure}[htbp]
\begin{center}
\includegraphics[scale=0.6]{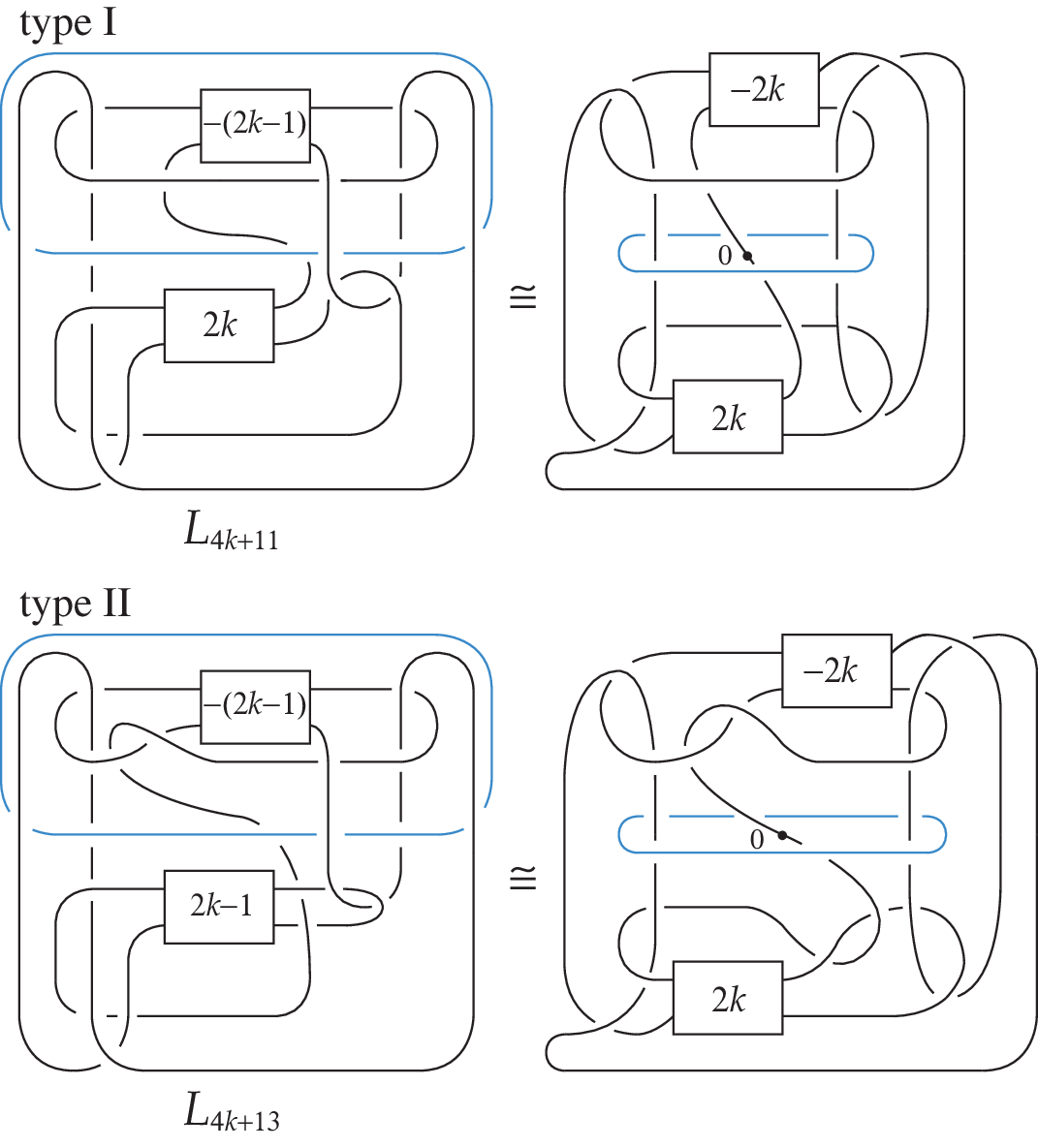} 
\label{ex}
\caption{prime component-preservingly amphicheiral link $L_n$}
\end{center}
\end{figure}

\noindent
{\bf Proof of Theorem \ref{thm:main}}\ 
By the righthand side of Figure 10, 
$L_n$ is a component-preservingly strongly $(-, +)$-amphicheiral link.

\medskip

The linking number of $L_n$, $\mathrm{lk}\ \! (L_n)$, 
is $3$ by a suitable orientation.
Let $c(\cdot)$ denote the minimal crossing number of a link.
Since
$$c(L_n)\ge c(S_n)+c(U)+2|\mathrm{lk}\ \! (L_n)|=n+6,$$
and the lefthand side of Figure 10 realizes the lower bound,
we have $c(L_n)=n+6$ and it is odd.

\medskip

Finally we show that $L_n$ is prime by using the Kauffman bracket.
Suppose that $L_n$ is not prime.
Then $L_n$ is a connected sum of two links such that
one is a Stoimenow knot $S_n$
and the other is a $2$-component link
with unknotted components and with linking number $3$
by Theorem \ref{thm:Stoimenow}.
Hence $\langle L_n\rangle$ should be divisible by $\langle S_n\rangle$
by Lemma \ref{lm:Kauffman} (4).
We compute $\langle L_n\rangle (\zeta)$ and
$\langle S_n\rangle (\zeta)$, where $\zeta$ is a primitive
$8$-th root of unity.
By Lemma \ref{lm:Kauffman} (4) and (8),
$|\langle L_n\rangle (\zeta)|$ should be divisible 
by $|\langle S_n\rangle (\zeta)|$.

\medskip

To compute $\langle S_n\rangle$ and $\langle L_n\rangle$,
we set $K=S_n$ and $L=L_n$, and
we denote the results of splicings by
$K_{00}$, $K_{0\infty}$, $K_{\infty 0}$, $K_{\infty \infty}$, 
$L_{00}$, $L_{0\infty}$, $L_{\infty 0}$ and $L_{\infty \infty}$, 
respectively as in Figure 11.
Here we drew only the type I case.
We can obtain the type II case in a similar way.

\begin{figure}[htbp]
\begin{center}
\includegraphics[scale=0.6]{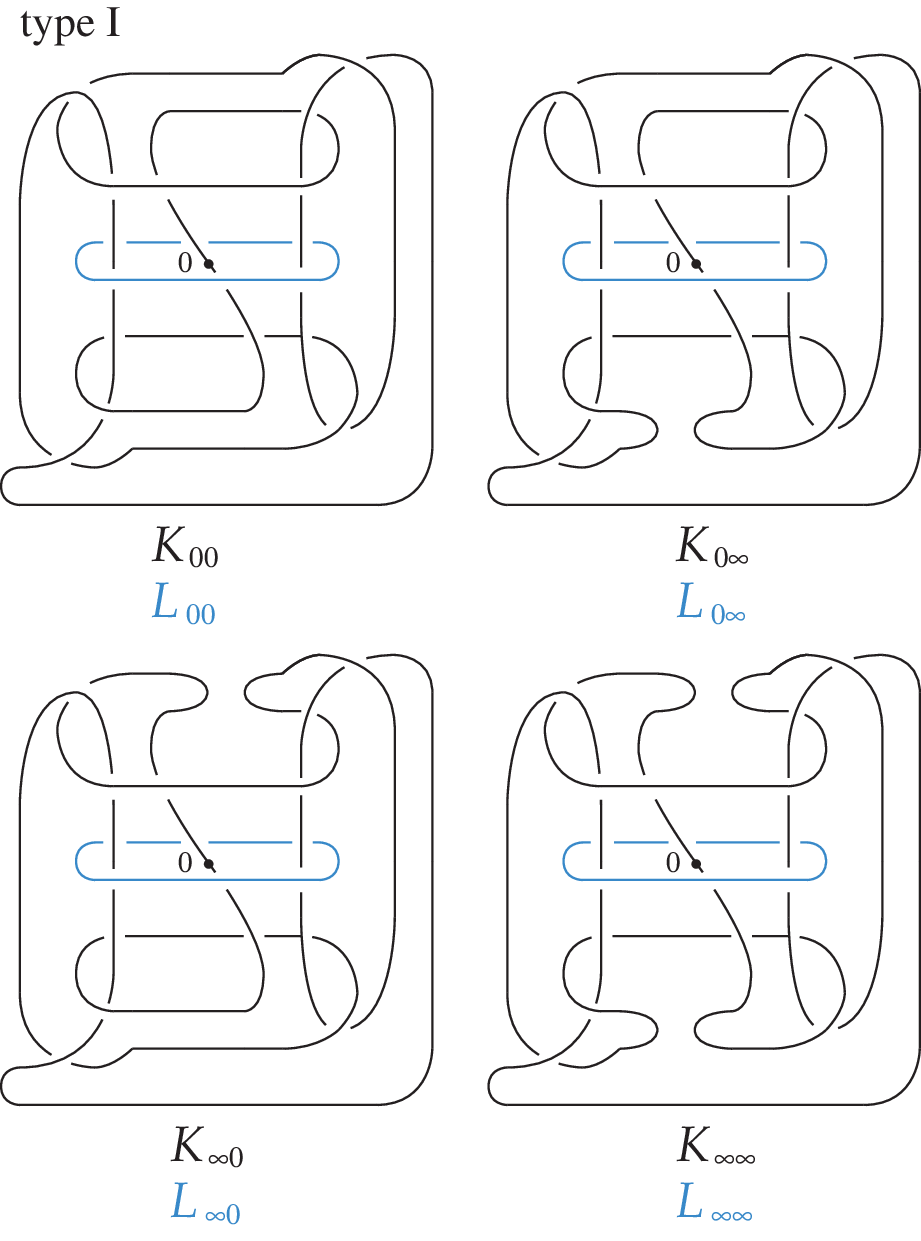} 
\label{splice2}
\caption{splices of $L_n$}
\end{center}
\end{figure}

\newpage

Then by Lemma \ref{lm:Tm} (1), we have:
\begin{equation}\label{eq:K}
\begin{matrix}
\langle K\rangle & = & \langle K_{00}\rangle
+A^{-2k}\alpha_{2k}(A)\langle K_{0\infty}\rangle
+A^{2k}\alpha_{-2k}(A)\langle K_{\infty 0}\rangle 
\hfill
\medskip\\
& &
+\alpha_{2k}(A)\alpha_{-2k}(A)\langle K_{\infty \infty}\rangle
\hfill
\end{matrix}
\end{equation}
and
\begin{equation}\label{eq:L}
\begin{matrix}
\langle L\rangle & = & \langle L_{00}\rangle
+A^{-2k}\alpha_{2k}(A)\langle L_{0\infty}\rangle
+A^{2k}\alpha_{-2k}(A)\langle L_{\infty 0}\rangle 
\hfill
\medskip\\
& &
+\alpha_{2k}(A)\alpha_{-2k}(A)\langle L_{\infty \infty}\rangle.
\hfill
\end{matrix}
\end{equation}
We can see that
$K_{00}$ and $K_{\infty \infty}$ are amphicheiral knot diagrams with writhe $0$,
$K_{0\infty}=(K_{\infty 0})^*$,
the writhe of $K_{0\infty}$ is $-10$,
the writhe of $K_{\infty 0}$ is $10$,
$L_{00}$ and $L_{\infty \infty}$
are $2$-component amphicheiral link diagrams with writhe $6$, 
$L_{0\infty}=(L_{\infty 0})^*$,
the writhe of $L_{0\infty}$ is $-4$,
and the writhe of $L_{\infty 0}$ is $16$.
By Lemma \ref{lm:Kauffman} (6),
we have 
$$K_{00}(A)=K_{00}(A^{-1}),\ 
K_{\infty \infty}(A)=K_{\infty \infty}(A^{-1}),\ 
K_{\infty 0}(A)=K_{0\infty}(A^{-1}),$$
$$L_{00}(A)=L_{00}(A^{-1}),\ 
L_{\infty \infty}(A)=L_{\infty \infty}(A^{-1}),\quad
\mbox{and}\quad 
L_{\infty 0}(A)=L_{0\infty}(A^{-1}).$$
By Lemma \ref{lm:Tm} (2),
$A^{2k}\alpha_{-2k}(A)$ can be obtained by replacing 
$A$ with $A^{-1}$ in $A^{-2k}\alpha_{2k}(A)$.
By straight calculations using 
Lemma \ref{lm:Kauffman} and Lemma \ref{lm:Tm},
we have:

\medskip

\noindent
(type I)
\begin{equation}\label{eq:typeIK}
\begin{matrix}
\langle K_{00}\rangle & = & 
A^{16}-4A^{12}+6A^8-7A^4+9-7A^{-4}+6A^{-8}-4A^{-12}+A^{-16},
\hfill
\medskip\\
\langle K_{0\infty}\rangle & = & 
-A^{18}+3A^{14}-5A^{10}+6A^6-7A^2+6A^{-2}
-5A^{-6}+4A^{-10}
\hfill \\
& & -A^{-14}+A^{-18},
\hfill
\medskip\\
\langle K_{\infty \infty}\rangle & = & 
A^{16}-3A^{12}+5A^8-6A^4+7-6A^{-4}+5A^{-8}-3A^{-12}+A^{-16}.
\hfill
\end{matrix}
\end{equation}
\begin{equation}\label{eq:typeIL}
\begin{matrix}
\langle L_{00}\rangle & = & 
-A^{20}+4A^{16}-8A^{12}+12A^8-16A^4+16
-16A^{-4}+12A^{-8}
\hfill \\
& & -8A^{-12}+4A^{-16}-A^{-20},
\hfill
\medskip\\
\langle L_{0\infty}\rangle & = & 
A^{22}-3A^{18}+6A^{14}-9A^{10}+12A^6-12A^2
+11A^{-2}-9A^{-6}
\hfill \\
& & +5A^{-10}-3A^{-14}-A^{-26},
\hfill
\medskip\\
\langle L_{\infty \infty}\rangle & = & 
-A^{20}+3A^{16}-7A^{12}+10A^8-13A^4+14
-13A^{-4}+10A^{-8}
\hfill \\
& & -7A^{-12}+3A^{-16}-A^{-20}.
\hfill
\end{matrix}
\end{equation}

\noindent
(type II)
\begin{equation}\label{eq:typeIIK}
\begin{matrix}
\langle K_{00}\rangle & = & 
-A^{20}+4A^{16}-9A^{12}+14A^8-17A^4+19
-17A^{-4}+14A^{-8}
\hfill \\
& & -9A^{-12}+4A^{-16}-A^{-20},
\hfill
\medskip\\
\langle K_{0\infty}\rangle & = & 
A^{22}-4A^{18}+10A^{14}-15A^{10}+19A^6-22A^2
+20A^{-2}-18A^{-6}
\hfill \\
& & +12A^{-10}-7A^{-14}+3A^{-18}-A^{-22},
\hfill
\medskip\\
\langle K_{\infty \infty}\rangle & = & 
-2A^{20}+6A^{16}-13A^{12}+21A^8-24A^4+28
-24A^{-4}+21A^{-8}
\hfill \\
& & -13A^{-12}+6A^{-16}-2A^{-20}.
\hfill
\end{matrix}
\end{equation}
\begin{equation}\label{eq:typeIIL}
\begin{matrix}
\langle L_{00}\rangle & = & 
A^{24}-5A^{20}+13A^{16}-24A^{12}+35A^8-44A^4
+46-44A^{-4}
\hfill \\
& & +35A^{-8}-24A^{-12}+13A^{-16}-5A^{-20}+A^{-24},
\hfill
\medskip\\
\langle L_{0\infty}\rangle & = & 
-A^{26}+4A^{22}-11A^{18}+20A^{14}-31A^{10}+40A^6
-42A^2+42A^{-2}
\hfill \\
& & -33A^{-6}+24A^{-10}-13A^{-14}+5A^{-18}-A^{-26}+A^{-30},
\hfill
\medskip\\
\langle L_{\infty \infty}\rangle & = & 
A^{24}-5A^{20}+14A^{16}-27A^{12}+38A^8-50A^4
+50-50A^{-4}
\hfill \\
& & +38A^{-8}-27A^{-12}+14A^{-16}-5A^{-20}+A^{-24}.
\hfill
\end{matrix}
\end{equation}

\medskip

We substitute $A=\zeta$ to (\ref{eq:K}) and (\ref{eq:L}).
We set $\zeta^2=\sqrt{-1}$.
By Lemma \ref{lm:Tm} (2) and the arguments above,
we have
\begin{equation}\label{eq:Kzeta}
\begin{matrix}
\langle K\rangle (\zeta) & = & \langle K_{00}\rangle (\zeta)
-4k\sqrt{-1}\langle K_{0\infty}\rangle (\zeta)
+4k^2\langle K_{\infty \infty}\rangle (\zeta)
\hfill
\end{matrix}
\end{equation}
and
\begin{equation}\label{eq:Lzeta}
\begin{matrix}
\langle L\rangle (\zeta) & = & \langle L_{00}\rangle (\zeta)
-4k\sqrt{-1}\langle L_{0\infty}\rangle (\zeta)
+4k^2\langle L_{\infty \infty}\rangle (\zeta).
\hfill
\end{matrix}
\end{equation}
By (\ref{eq:typeIK}), (\ref{eq:typeIL}), (\ref{eq:typeIIK}) 
and (\ref{eq:typeIIL}), we have

\medskip

\noindent
(type I)
\begin{equation}\label{eq:typeIKzeta}
\begin{matrix}
\langle K_{00}\rangle (\zeta) & = & 45,
\hfill
\medskip\\
\langle K_{0\infty}\rangle (\zeta) & = & -39\sqrt{-1},
\hfill
\medskip\\
\langle K_{\infty \infty}\rangle (\zeta) & = & 37.
\hfill
\end{matrix}
\end{equation}
\begin{equation}\label{eq:typeILzeta}
\begin{matrix}
\langle L_{00}\rangle (\zeta) & = & 98,
\hfill
\medskip\\
\langle L_{0\infty}\rangle (\zeta) & = & -70\sqrt{-1},
\hfill
\medskip\\
\langle L_{\infty \infty}\rangle (\zeta) & = & 82.
\hfill
\end{matrix}
\end{equation}

\noindent
(type II)
\begin{equation}\label{eq:typeIIKzeta}
\begin{matrix}
\langle K_{00}\rangle (\zeta) & = & 109,
\hfill
\medskip\\
\langle K_{0\infty}\rangle (\zeta) & = & -132\sqrt{-1},
\hfill
\medskip\\
\langle K_{\infty \infty}\rangle (\zeta) & = & 160.
\hfill
\end{matrix}
\end{equation}
\begin{equation}\label{eq:typeIILzeta}
\begin{matrix}
\langle L_{00}\rangle (\zeta) & = & 290,
\hfill
\medskip\\
\langle L_{0\infty}\rangle (\zeta) & = & -264\sqrt{-1},
\hfill
\medskip\\
\langle L_{\infty \infty}\rangle (\zeta) & = & 320.
\hfill
\end{matrix}
\end{equation}

By (\ref{eq:Kzeta}), (\ref{eq:Lzeta}),
(\ref{eq:typeIKzeta}), (\ref{eq:typeILzeta}), (\ref{eq:typeIIKzeta}) 
and (\ref{eq:typeIILzeta}), we have

\medskip

\noindent
(type I)
\begin{equation*}\label{eq:typeI}
\begin{matrix}
\langle K\rangle (\zeta) & = & 148k^2-156k+45,
\hfill
\medskip\\
\langle L\rangle (\zeta) & = & 328k^2-280k+98.
\end{matrix}
\end{equation*}

\noindent
(type II)
\begin{equation*}\label{eq:typeII}
\begin{matrix}
\langle K\rangle (\zeta) & = & 640k^2-528k+109,
\hfill
\medskip\\
\langle L\rangle (\zeta) & = & 1280k^2-1056k+290.
\end{matrix}
\end{equation*}
Note that $148k^2-156k+45$ and $640k^2-528k+109$ are odd 
and $328k^2-280k+98$ and $1280k^2-1056k+290$ are 
of the form $2\times$(odd), and
they are positive for $k\ge 1$.
Hence if $148k^2-156k+45$ divides $328k^2-280k+98$
($640k^2-528k+109$ divides $1280k^2-1056k+290$, respectively),
then $148k^2-156k+45$ divides $164k^2-140k+49$
($640k^2-528k+109$ divides $640k^2-528k+145$, respectively),
and the quantity is odd.

\medskip

\noindent
(type I)

\medskip

Suppose that $164k^2-140k+49$ is divisible by $148k^2-156k+45$.
Since
$$(164k^2-140k+49)-(148k^2-156k+45)
=16k^2+16k+4>0,$$
the quantity is not $1$.
Since
$$3(148k^2-156k+45)-(164k^2-140k+49)
=280k^2-328k+86>0,$$
the quantity is not greater than $1$.
It is a contradiction.

\medskip

\noindent
(type II)

\medskip

Suppose that $640k^2-528k+145$ is divisible by $640k^2-528k+109$.
Since
$$(640k^2-528k+145)-(640k^2-528k+109)
=36>0,$$
the quantity is not $1$.
Since
$$3(640k^2-528k+109)-(640k^2-528k+145)
=1280k^2-1056k+182>0,$$
the quantity is not greater than $1$.
It is a contradiction.
\qed

\begin{re}\label{re:Remark}
{\rm
In \cite{Ko}, the second author computes
the J polynomials, which are modified Jones polynomials, 
of $S_n$ and $L_n$ explicitly.
The J polynomial is an invariant of unoriented links.}
\end{re}

\section{Non-invertibility of Stoimenow knots}\label{sec:noninv}

In this section, we show that a Stoimenow knot $S_n$ is not invertible
by using the Alexander polynomials.
Since $S_n$ is $(-)$-amphicheiral, 
we show that it is not $(+)$-amphicheiral, which is equivalent to that it is not invertble.

\medskip

Let $L$ be a link, and ${\mit \Delta}_L(t)\in \mathbb{Z}[t, t^{-1}]$ 
the Alexander polynomial of $L$.
For two elements $A$ and $B$ in 
$\mathbb{Z}[t, t^{-1}]$
($(\mathbb{Z}/d\mathbb{Z})[t, t^{-1}]$, 
respectively), we denote by 
$A\doteq B$ ($A\doteq_d B$, respectively)
if they are equal up to multiplications of trivial units.
A one variable Laurent polynomial $r(t)\in \mathbb{Z}[t^{\pm 1}]$
is {\it of type $X$} if there are integers
$n\ge 0$ and $\lambda \ge 3$ such that $\lambda$ is odd, and
$f_i(t)\in \mathbb{Z}[t, t^{-1}]$
$(i=0, 1, \ldots, n)$
such that
$f_i(t)\doteq f_i(t^{-1})$, 
$|f_i(1)|=1$, and for $i>0$,
$f_i(t)\doteq_2 f_0(t)^{2^i}p_{\lambda}(t)^{2^{i-1}}$
where $p_{\lambda}(t)=(t^{\lambda}-1)/(t-1)$, and
\begin{equation}\label{eq:r}
r(t)\doteq 
\left\{
\begin{array}{ll}
f_0(t)^2 & (n=0),
\medskip\\
f_0(t)^2f_1(t)\cdots f_n(t) & (n\ge 1).
\end{array}
\right.
\end{equation}
R.~Hartley \cite{Ha},
R.~Hartley and A.~Kawauchi \cite{HK}, and A.~Kawauchi \cite{Kw1}
gave necessary conditions on the Alexander polynomials
of amphicheiral knots.

\begin{lm}\label{lm:HK}
{\rm (Hartley \cite{Ha}; Hartley and Kawauchi \cite{HK}; Kawauchi \cite{Kw1})}
\begin{enumerate}
\item[(1)]
Let $K$ be a $(-)$-amphicheiral knot.
Then there exists an element $f(t)\in \mathbb{Z}[t, t^{-1}]$
such that
$|f(1)|=1$, $f(t^{-1})\doteq f(-t)$, and
$${\mit \Delta}_K(t^2)\doteq f(t)f(t^{-1}).$$

\item[(2)]
Let $K$ be a $(+)$-amphicheiral knot.
Then there exist
$r_j(t)\in \mathbb{Z}[t, t^{-1}]$ of type $X$
and a positive odd number $\alpha_j$\ $(j=1, \ldots, m)$
such that
$${\mit \Delta}_K(t)\doteq \prod_{j=1}^m r_j(t^{\alpha_j}).$$
In particular, if $K$ is hyperbolic, then
we can take $m=1$ and $\alpha_1=1$.

\end{enumerate}
\end{lm}

We generalize Stoimenow knots as in Figure 12.
The lefthand side is called a {\it generalized Stoimenow link of type I},
and is denoted by $S^1_{p,q}$.
The righthand side is called a {\it generalized Stoimenow link of type II},
and is denoted by $S^2_{r,s}$.
The numbers in rectangles are the numbers of half twists.
We note that $S^1_{2k,2k}=S_{4k+11}$ and $S^2_{k,k}=S_{4k+13}$.
We denote the Alexander polynomials (the Conway polynomials) of 
$S^1_{p,q}$ and $S^2_{r,s}$ by 
${\mit \Delta}^{(1)}_{p,q}(t)$ and ${\mit \Delta}^{(2)}_{r,s}(t)$
($\nabla^{(1)}_{p,q}(z)$ and $\nabla^{(2)}_{r,s}(z)$), respectively.
We compute ${\mit \Delta}^{(1)}_{2k,2k}(t)$ and ${\mit \Delta}^{(2)}_{k,k}(t)$
as the mod $2$ Alexander polynomials.

\begin{figure}[htbp]
\begin{center}
\includegraphics[scale=0.6]{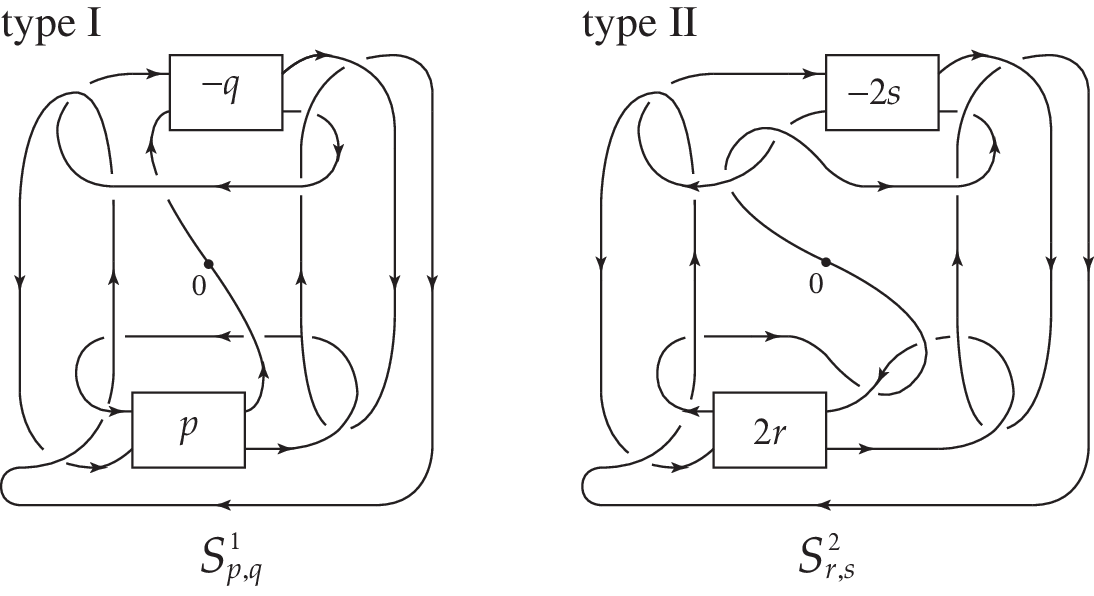} 
\label{general}
\caption{generalized Stoimenow links $S_{p,q}^1$ and $S_{r,s}^2$}
\end{center}
\end{figure}

\newpage

\begin{lm}\label{lm:StAlex}
The Alexander and the mod $2$ Alexander polynomials of 
$S^1_{2k,2k}$ and $S^2_{k,k}$ are as follows : 
\begin{eqnarray*}
(t+1)^2 {\mit \Delta}^{(1)}_{2k,2k}(t)
& \doteq_2 & t^{4k+6}+t^{4k+5}+t^{4k+4}+t^{4k+2}+t^{4k-1}+t^7+t^4+t^2+t+1
\bigskip\\
& =_2 & (t^2+t+1)^2(t^{4k+2}+t^{4k+1}+t^{4k-1}+t^3+t+1).
\bigskip\\
{\mit \Delta}^{(2)}_{k,k}(t)
& \doteq & 
t^3(-t^6+9t^5-26t^4+37t^3-26t^2+9t-1)\\
&  & 
-2kt^2(t-1)^2(2t^6-7t^5+15t^4-18t^3+15t^2-7t+2)\\
&  & 
+k^2(t-1)^2(t^{10}-3t^9+7t^8-17t^7+32t^6-40t^5+32t^4-17t^3\\
&  & 
+7t^2-3t+1)
\bigskip\\
& \doteq_2 & 
\left\{
\begin{array}{cl}
t^6+t^5+t^3+t+1 & (\mbox{$k$ is even}),
\medskip\\
t^{12}+t^{11}+t^9+t^7+t^6+t^5+t^3+t+1 & (\mbox{$k$ is odd}).
\end{array}
\right.
\end{eqnarray*}

\end{lm}

\noindent
{\bf Proof}\ \ We have the following relations on the Conway polynomials 
from the skein relation in Subsection \ref{ssec:Alexander} : 

\begin{equation}\label{eq:rel1}
\left\{
\begin{array}{ccccl}
\nabla_{p,q}^{(1)}(z) & - & \nabla_{p-2,q}^{(1)}(z) & = & z\nabla_{p-1,q}^{(1)}(z)
\medskip \\
\nabla_{p,q-2}^{(1)}(z) & - & \nabla_{p,q}^{(1)}(z) & = & z\nabla_{p,q-1}^{(1)}(z),
\end{array}
\right.
\end{equation}
and
\begin{equation}\label{eq:rel2}
\left\{
\begin{array}{ccccl}
\nabla_{r-1,s}^{(2)}(z) & - & \nabla_{r,s}^{(2)}(z) & = & z\nabla_{\infty,s}^{(2)}(z)
\medskip \\
\nabla_{r,s}^{(2)}(z) & - & \nabla_{r,s-1}^{(2)}(z) & = & z\nabla_{r,\infty}^{(2)}(z).
\end{array}
\right.
\end{equation}
For the meaning of $\infty$, see Figure 6.

\medskip

\noindent
(type I)

\medskip

From (\ref{eq:rel1}), we have : 

\begin{equation}\label{eq:rel3}
\left\{
\begin{array}{ccl}
{\mit \Delta}_{p,q}^{(1)}(t)-t^{\frac 12}{\mit \Delta}_{p-1,q}^{(1)}(t) & 
= & (-t^{-\frac 12})^{p-1}({\mit \Delta}_{1,q}^{(1)}(t)-t^{\frac 12}{\mit \Delta}_{0,q}^{(1)}(t))
\medskip \\
{\mit \Delta}_{p,q}^{(1)}(t)+t^{-\frac 12}{\mit \Delta}_{p-1,q}^{(1)}(t) & 
= & (t^{-\frac 12})^{p-1}({\mit \Delta}_{1,q}^{(1)}(t)+t^{-\frac 12}{\mit \Delta}_{0,q}^{(1)}(t)),
\end{array}
\right.
\end{equation}
and 
\begin{equation}\label{eq:rel4}
\left\{
\begin{array}{ccl}
{\mit \Delta}_{p,q}^{(1)}(t)+t^{\frac 12}{\mit \Delta}_{p,q-1}^{(1)}(t) & 
= & (t^{-\frac 12})^{q-1}({\mit \Delta}_{p,1}^{(1)}(t)+t^{\frac 12}{\mit \Delta}_{p,0}^{(1)}(t))
\medskip \\
{\mit \Delta}_{p,q}^{(1)}(t)-t^{-\frac 12}{\mit \Delta}_{p,q-1}^{(1)}(t) & 
= & (-t^{-\frac 12})^{q-1}({\mit \Delta}_{p,1}^{(1)}(t)-t^{-\frac 12}{\mit \Delta}_{p,0}^{(1)}(t)).
\end{array}
\right.
\end{equation}

\medskip

From (\ref{eq:rel3}) and (\ref{eq:rel4}), we have : 

\begin{equation}\label{eq:rel5}
\left\{
\begin{array}{ccl}
(t^{\frac 12}+t^{-\frac 12}){\mit \Delta}_{p,q}^{(1)}(t) & 
= & (t^{\frac p2}-(-1)^pt^{-\frac p2}){\mit \Delta}_{1,q}^{(1)}(t)
\medskip \\
& & +(t^{\frac{p-1}{2}}+(-1)^pt^{-\frac{p-1}{2}}){\mit \Delta}_{0,q}^{(1)}(t)
\medskip \\
-(t^{\frac 12}+t^{-\frac 12}){\mit \Delta}_{p,q}^{(1)}(t) & 
= & ((-1)^qt^{\frac q2}-t^{-\frac q2}){\mit \Delta}_{p,1}^{(1)}(t)
\medskip \\
& & -((-1)^qt^{\frac{q-1}{2}}+t^{-\frac{q-1}{2}}){\mit \Delta}_{p,0}^{(1)}(t).
\end{array}
\right.
\end{equation}

\medskip

From (\ref{eq:rel5}), if $p=q=2k$, then we have a skein relation 
among the Alexander polynomials of 
$S^1_{2k,2k}$, $S^1_{0,0}$, $S^1_{1,0}$, $S^1_{0,1}$ and $S^1_{1,1}$
(cf.\ Figure 13) : 

\begin{equation}\label{eq:type11}
\begin{array}{ccl}
\left( t^{\frac 12}+t^{-\frac 12}\right)^2 {\mit \Delta}^{(1)}_{2k,2k}(t)
& = & \left( t^{k-\frac 12}+t^{-k+\frac 12}\right)^2 {\mit \Delta}^{(1)}_{0,0}(t)
\medskip\\
& & 
-( t^k+t^{-k})\left( t^{k-\frac 12}+t^{-k+\frac 12}\right) 
({\mit \Delta}^{(1)}_{1,0}(t)-{\mit \Delta}^{(1)}_{0,1}(t))
\medskip\\
& & -( t^k+t^{-k})^2 {\mit \Delta}^{(1)}_{1,1}(t).
\end{array}
\end{equation}

\medskip

Since $S^1_{1,0}$ and $S^1_{0,1}$ are $2$-component links 
with $S^1_{0,1}=-(S^1_{1,0})^*$, and (\ref{eq:type11}), 
we have $\nabla^{(1)}_{0,1}(z)=-\nabla^{(1)}_{1,0}(z)$
and ${\mit \Delta}^{(1)}_{0,1}(t)=-{\mit \Delta}^{(1)}_{1,0}(t)$
by Lemma \ref{lm:Conway}, and

\begin{equation}\label{eq:type12}
\begin{array}{ccl}
\left( t^{\frac 12}+t^{-\frac 12}\right)^2 {\mit \Delta}^{(1)}_{2k,2k}(t)
& = & \left( t^{k-\frac 12}+t^{-k+\frac 12}\right)^2 {\mit \Delta}^{(1)}_{0,0}(t)
\medskip\\
& & 
-2( t^k+t^{-k})\left( t^{k-\frac 12}+t^{-k+\frac 12}\right) 
{\mit \Delta}^{(1)}_{1,0}(t)
\medskip\\
& & 
-( t^k+t^{-k})^2 {\mit \Delta}^{(1)}_{1,1}(t).
\end{array}
\end{equation}

\medskip

Since $S^1_{0,0}=8_{18}$,
$$\begin{array}{rll}
{\mit \Delta}^{(1)}_{0,0}(t)={\mit \Delta}_{8_{18}}(t)
& = & -t^3+5t^2-10t+13-10t^{-1}+5t^{-2}-t^{-3}\\
& =_2 & t^3+t^2+1+t^{-2}+t^{-3},
\medskip \\
{\mit \Delta}^{(1)}_{1,1}(t)
& = & -t^{-3}(t^3-1)^2=_2t^3+t^{-3},
\end{array}$$
and (\ref{eq:type12}), we have
$$\begin{array}{rcl}
(t+1)^2 {\mit \Delta}^{(1)}_{2k,2k}(t)
& \doteq_2 & t^{4k+6}+t^{4k+5}+t^{4k+4}+t^{4k+2}+t^{4k-1}+t^7+t^4+t^2+t+1
\medskip\\
& =_2 & (t^2+t+1)^2(t^{4k+2}+t^{4k+1}+t^{4k-1}+t^3+t+1).
\end{array}$$

\noindent
(type II)

\medskip

From (\ref{eq:rel2}), we have : 

\begin{equation}\label{eq:rel6}
\left\{
\begin{array}{ccl}
\nabla_{r,s}^{(2)}(z) & 
= & \nabla_{0,s}^{(2)}(z)-rz\nabla_{\infty,s}^{(2)}(z)
\medskip \\
\nabla_{r,s}^{(2)}(z) & 
= & \nabla_{r,0}^{(2)}(z)+sz\nabla_{r,\infty}^{(2)}(z).
\end{array}
\right.
\end{equation}

\medskip

From (\ref{eq:rel6}), we have : 

\begin{equation*}
\nabla_{r,s}^{(2)}(z)=
\nabla_{0,0}^{(2)}(z)-rz\nabla_{\infty,0}^{(2)}(z)+sz\nabla_{0,\infty}^{(2)}(z)
-rsz^2\nabla_{\infty,\infty}^{(2)}(z).
\end{equation*}
In particular, if $r=s=k$, then we have a skein relation 
among the Conway polynomials of 
$S^2_{k,k}$, $S^2_{0,0}$, $S^2_{0,\infty}$, $S^2_{\infty,0}$ and $S^2_{\infty,\infty}$
(cf.\ Figure 14) : 
\begin{equation}\label{eq:rel7}
\nabla_{k,k}^{(2)}(z)=
\nabla_{0,0}^{(2)}(z)+kz(\nabla_{0,\infty}^{(2)}(z)-\nabla_{\infty,0}^{(2)}(z))
-k^2z^2\nabla_{\infty,\infty}^{(2)}(z).
\end{equation}

\medskip

Since $S^2_{0,\infty}$ and $S^2_{\infty,0}$ are $2$-component links 
with $S^2_{0,\infty}=-(S^2_{\infty,0})^*$, 
$$\begin{array}{cll}
\nabla^{(2)}_{0,0}(z)
& = & -z^6+3z^4+z^2+1,
\medskip \\
\nabla^{(2)}_{0,\infty}(z)
& = & -2z^7-5z^5-5z^3-2z,
\medskip \\
\nabla^{(2)}_{\infty,\infty}(z)
& = & -z^{10}-7z^8-18z^6-15z^4-4z^2,
\end{array}$$
and (\ref{eq:rel7}), 
we have $\nabla^{(2)}_{0,\infty}(z)=-\nabla^{(2)}_{\infty,0}(z)$
by Lemma \ref{lm:Conway}, and
\begin{eqnarray*}
{\mit \Delta}^{(2)}_{k,k}(t)
& \doteq & 
t^3(-t^6+9t^5-26t^4+37t^3-26t^2+9t-1)\\
&  & 
-2kt^2(t-1)^2(2t^6-7t^5+15t^4-18t^3+15t^2-7t+2)\\
&  & 
+k^2(t-1)^2(t^{10}-3t^9+7t^8-17t^7+32t^6-40t^5+32t^4-17t^3\\
&  & 
+7t^2-3t+1)
\bigskip\\
& \doteq_2 & 
\left\{
\begin{array}{cl}
t^6+t^5+t^3+t+1 & (\mbox{$k$ is even}),
\medskip\\
t^{12}+t^{11}+t^9+t^7+t^6+t^5+t^3+t+1 & (\mbox{$k$ is odd}).\ 
\qed
\end{array}
\right.
\end{eqnarray*}

\begin{figure}[htbp]
\begin{center}
\includegraphics[scale=0.6]{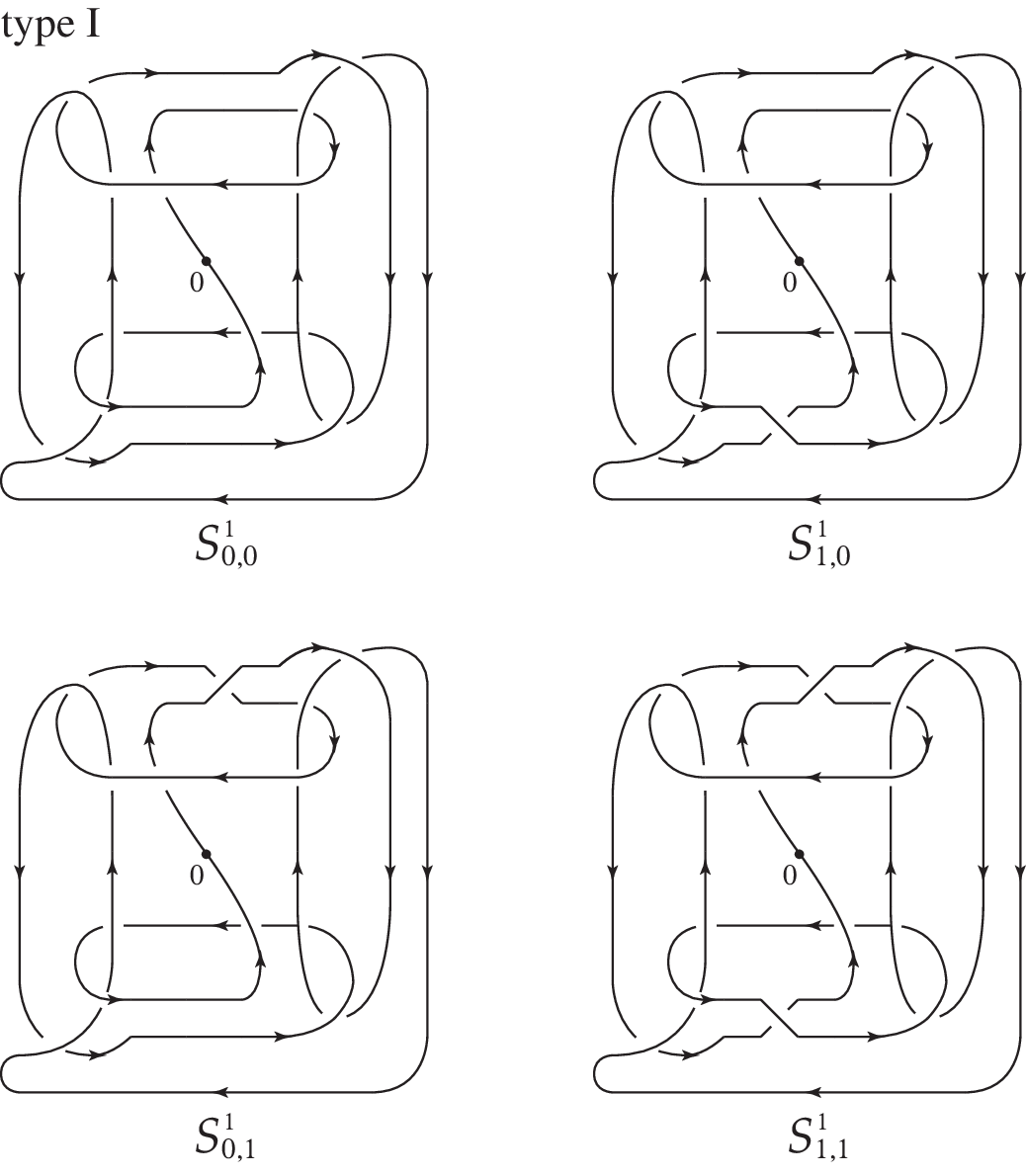} 
\label{type1}
\caption{$S_{0,0}^1$, $S_{1,0}^1$, $S_{0,1}^1$ and $S_{1,1}^1$}
\end{center}
\end{figure}

\newpage

\begin{figure}[htbp]
\begin{center}
\includegraphics[scale=0.6]{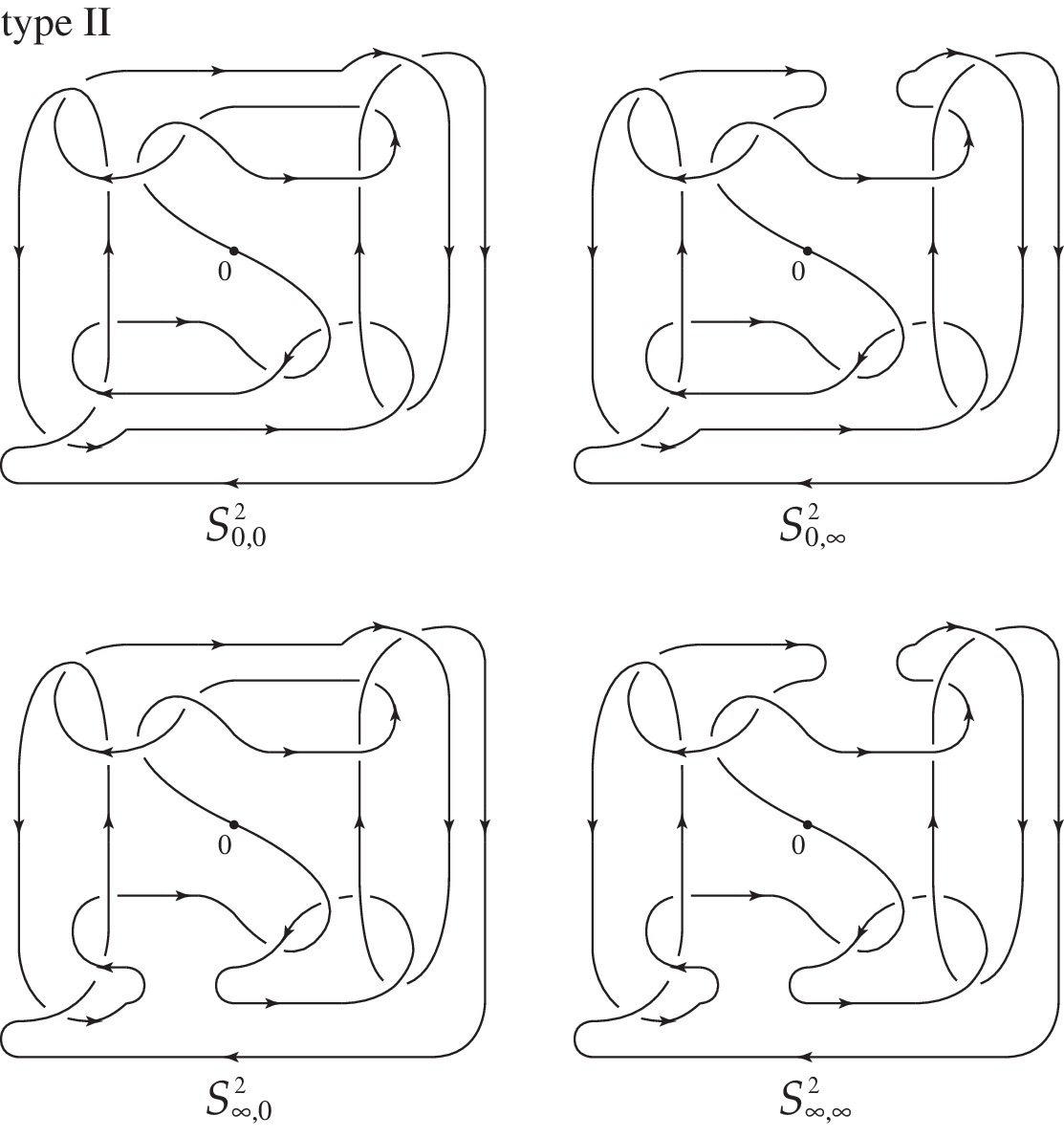} 
\label{type2}
\caption{$S_{0,0}^2$, $S_{0,\infty}^2$, $S_{\infty,0}^2$ and $S_{\infty,\infty}^2$}
\end{center}
\end{figure}

\newpage

Every element $f\in (\mathbb{Z}/2\mathbb{Z})[t, t^{-1}]$ is of the form : 
$$f=t^{k_d}+t^{k_{d-1}}+\cdots +t^{k_1}+t^{k_0}$$
where $k_0, \ldots, k_d$ are integers such that
$k_0<k_1<\cdots <k_{d-1}<k_d$.
Then we define the {\it mod $2$ trace}, 
denoted by $\mathrm{tr}_2(f)\in \mathbb{Z}/2\mathbb{Z}=\{0, 1\}$, as : 
$$\mathrm{tr}_2(f)
=\left\{
\begin{array}{cl}
1 & (k_d-k_{d-1}=1),\medskip\\
0 & (k_d-k_{d-1}\ge 2).
\end{array}
\right.$$
For $f_1, f_2\in (\mathbb{Z}/2\mathbb{Z})[t, t^{-1}]$, 
$\mathrm{tr}_2(f_1f_2)=\mathrm{tr}_2(f_1)+\mathrm{tr}_2(f_2)$.
There exists an element $g\in (\mathbb{Z}/2\mathbb{Z})[t, t^{-1}]$
such that $f=g^2$ if and only if 
every $k_i$\ $(i=0, \ldots, d)$ is even.
Then we call $f$ a {\it square polynomial}, and 
we have 
$$g=t^{k_d/2}+\cdots +t^{k_1/2}+t^{k_0/2}$$ 
and $\mathrm{tr}_2(f)=0$.

\begin{lm}\label{lm:square}
Let $r(t)$ be of type X as in (\ref{eq:r}), and $\alpha$ a positive odd integer.
\begin{enumerate}
\item[(1)]
If $n=0$, then $r(t^{\alpha})$ is a square polynomial.
If $n\ge 1$, then $r(t^{\alpha})$ is of the form : 
$$r(t^{\alpha})=g^2p_{\lambda}(t^{\alpha})$$
where $g\in (\mathbb{Z}/2\mathbb{Z})[t, t^{-1}]$
and $p_{\lambda}(t)=(t^{\lambda}-1)/(t-1)$.

\item[(2)]
$\mathrm{tr}_2(r(t^{\alpha}))=1$ if and only if 
$n\ge 1$ and $\alpha=1$.

\end{enumerate}
\end{lm}

Let $\zeta_m$ be a primitive $m$-th root of unity, and
$\mathbf{\Phi}_m(t)\in \mathbb{Z}[t]$ the $m$-th cyclotomic polynomial
defined by
$$\mathbf{\Phi}_m(t)=\prod_{
{\scriptstyle 1\le i\le m-1}
\atop
{\scriptstyle \gcd(i, m)=1}}(t-\zeta_m^i).$$
The cyclotomic polynomial is a monic symmetric irreducible polynomial over $\mathbb{Z}$.
For a prime $q$ and a positive integer $r$, 
$$\mathbf{\Phi}_{q^r}(t)=\frac{t^{q^r}-1}{t^{q^{r-1}}-1}
=t^{q^{r-1}(q-1)}+t^{q^{r-1}(q-2)}+\cdots +t^{q^{r-1}}+1.$$
Since
$$t^m-1=\prod_{d\ge 1, d|m}\mathbf{\Phi}_d(t),$$
we have
\begin{equation}\label{eq:cyclotomic}
p_{\lambda}(t^{\alpha})
=\frac{t^{\alpha \lambda}-1}{t^{\alpha}-1}
=\prod_{d|\alpha \lambda, d\not\ \! | \alpha}\mathbf{\Phi}_d(t).
\end{equation}

\begin{theo}\label{th:Stoimenow}
A Stoimenow knot $S_n$ is not invertible.
\end{theo}

\noindent
{\bf Proof}\ \ 
We show that both $S^1_{2k,2k}$ and $S^2_{k,k}$ with $k\ge 1$
are not $(+)$-amphicheiral.

\medskip

\noindent
(type I)

\medskip

Suppose that ${\mit \Delta}^{(1)}_{2k,2k}(t)$ satisfies
the condition in Lemma \ref{lm:HK} (2).

\medskip

We set 
$$h=_2 
t^{4k+2}+t^{4k+1}+t^{4k-1}+t^3+t+1,$$
and $m=q^r$ with an odd prime $q\ge 3$ and $r\ge 1$.
Then 
\begin{equation}\label{eq:h}
(t+1)^2{\mit \Delta}^{(1)}_{2k,2k}(t)\doteq_2 (t^2+t+1)^2h.
\end{equation}

\noindent
{\bf Claim 1}\ 
{\it $\mathbf{\Phi}_m(t)$ is a mod $2$ divisor of $h$ 
only if $m=3, 5$ or $9$.}

\medskip

\noindent
{\bf Proof}\ \ 
Take $Q(t), R(t)\in (\mathbb{Z}/2\mathbb{Z})[t, t^{-1}]$ 
such that $h=_2\mathbf{\Phi}_m(t)Q(t)+R(t)$.
We can take $R(t)$ of the form : 
$$R(t)=_2t^{d+3}+t^{d+2}+t^d+t^3+t+1$$
where $-m/2<d<m/2$.
The span of $R(t)$ is less than $m/2+3$.

\medskip

\noindent
{\bf Case 1}\ $r\ge 2$ except the case $(q, r)=(3, 2)$.

\medskip

Since the degree of $\mathbf{\Phi}_m(t)$ is $q^{r-1}(q-1)$
which is greater than $q^r/2+3$, $R(t)=0$ should be hold.
However it does not occur.

\bigskip

\noindent
{\bf Case 2}\ $(q, r)=(3, 2)$ ($m=9$).

\medskip

$R(t)$ is not mod $2$ divisible by $\mathbf{\Phi}_9(t)=t^6+t^3+1$ 
except the case $d=4$.

\bigskip

\noindent
{\bf Case 3}\ $r=1$.

\medskip

We check only the cases $m=3, 5$ and $7$.
The case $m=7$ does not occur.
Hence we have the result.
\qed

\bigskip

\noindent
{\bf Claim 2}\ 
{\it
$h$ is mod $2$ divisible by $\mathbf{\Phi}_3(t)$ if and only if $k\equiv 0\ (\mathrm{mod}\ \! 3)$.
$h$ is mod $2$ divisible by $\mathbf{\Phi}_5(t)$ if and only if $k\equiv 1\ (\mathrm{mod}\ \! 5)$.
$h$ is mod $2$ divisible by $\mathbf{\Phi}_9(t)$ if and only if $k\equiv -1\ (\mathrm{mod}\ \! 9)$.}

\bigskip

\noindent
{\bf Proof}\ \ 
$h$ is mod $2$ divisible by $\mathbf{\Phi}_3(t)$ if and only if 
$4k+1\equiv 1\ (\mathrm{mod}\ \! 3)$ 
which is equivalent to $k\equiv 0\ (\mathrm{mod}\ \! 3)$.

\medskip

$h$ is mod $2$ divisible by $\mathbf{\Phi}_5(t)$ if and only if 
$4k+1\equiv 0\ (\mathrm{mod}\ \! 5)$ and $4k-1\equiv 3\ (\mathrm{mod}\ \! 5)$
which is equivalent to $k\equiv 1\ (\mathrm{mod}\ \! 5)$.

\medskip

$h$ is mod $2$ divisible by $\mathbf{\Phi}_9(t)$ if and only if 
$4k+1\equiv 6\ (\mathrm{mod}\ \! 9)$
which is equivalent to $k\equiv -1\ (\mathrm{mod}\ \! 9)$.
\qed

\bigskip

\noindent
{\bf Claim 3}\ 
{\it $\mathbf{\Phi}_{15}(t)$ is a mod $2$ divisor of $h$
if and only if $k\equiv -5\ (\mathrm{mod}\ \! 15)$.
$\mathbf{\Phi}_{45}(t)$ is not a mod $2$ divisor of $h$.
}

\bigskip

\noindent
{\bf Proof}\ \ 
For $\mathbf{\Phi}_{15}(t)=t^8-t^7+t^5-t^4+t^3-t+1$, 
we only check the cases $d=\pm 5, \pm 6$ and $\pm 7$.
For the cases, $R(t)$ is mod $2$ divisible by $\mathbf{\Phi}_{15}(t)$
if and only if $4k-1\equiv -6\ (\mathrm{mod}\ \! 15)$
which is equivalent to $k\equiv -5\ (\mathrm{mod}\ \! 15)$.

\medskip

For $\mathbf{\Phi}_{45}(t)=t^{24}-t^{21}+t^{15}-t^{12}+t^9-t^3+1$,
we only check the cases $d=\pm 21$ and $\pm 22$.
For the cases, $R(t)$ is not mod $2$ divisible by $\mathbf{\Phi}_{45}(t)$.
\qed

\bigskip

\noindent
{\bf Claim 4}\ 
{\it
$p_{\lambda}(t^{\alpha})$ is a mod $2$ divisor of $h$
only if $p_3(t)=\mathbf{\Phi}_3(t)=t^2+t+1$, 
$p_5(t)=\mathbf{\Phi}_5(t)=t^4+t^3+t^2+t+1$ or $p_3(t^3)=\mathbf{\Phi}_9(t)=t^6+t^3+1$.}

\bigskip

\noindent
{\bf Proof}\ \ 
By Claim 1, Claim 2, Claim 3 and (\ref{eq:cyclotomic}), we have the result.
\qed

\bigskip

By Lemma \ref{lm:StAlex}, we have $\mathrm{tr}_2({\mit \Delta}^{(1)}_{2k,2k}(t))=1$.
By Lemma \ref{lm:square}, Claim 1, Claim 2, Claim 3, Claim 4
and (\ref{eq:h}), $h$
is of the form : 
$$h\doteq_2 g^2p_3(t),\ g^2p_5(t)\ \mbox{or}\ g^2p_5(t)p_3(t^3)$$
for some $g\in (\mathbb{Z}/2\mathbb{Z})[t, t^{-1}]$.
However we have
$$\frac{h}{t^2+t+1}=_2t^{4k}+\cdots+t^5+t^4+t^2+1$$
for $k\equiv 0\ (\mathrm{mod}\ \! 3), k\ge 3$,
$$\frac{h}{t^4+t^3+t^2+t+1}=_2t^{4k-2}+\cdots+t^3+t^2+1$$
for $k\equiv 1\ (\mathrm{mod}\ \! 5), k\ge 6$,
and
$$\frac{h}{(t^4+t^3+t^2+t+1)(t^6+t^3+1)}=_2
t^{4k-8}+\cdots+t^5+t^2+1$$
for $k\equiv 26\ (\mathrm{mod}\ \! 45), k\ge 26$
are not square polynomials.
It is a contradiction.

\bigskip

\noindent
(type II)

\medskip

Suppose that ${\mit \Delta}^{(2)}_{k,k}(t)$ satisfies
the condition in Lemma \ref{lm:HK} (2).

\medskip

By Lemma \ref{lm:StAlex}, we have $\mathrm{tr}_2({\mit \Delta}^{(2)}_{k,k}(t))=1$.
By Lemma \ref{lm:square}, there exists an odd $\lambda \ge 3$ such that 
$p_{\lambda}(t)$ is a mod $2$ divisor of ${\mit \Delta}^{(2)}_{k,k}(t)$.
If $k$ is odd, then there is no such $\lambda$
(Check only the cases $\lambda=3, 5, 7, 9, 11$).
Hence we suppose that $k$ is even.
Since
$${\mit \Delta}^{(2)}_{k,k}(t)\doteq_2
(t^2+t+1)^3,$$
we have $\lambda=3$.
By the forms (\ref{eq:r}) and Lemma \ref{lm:HK} (2),
${\mit \Delta}^{(2)}_{k,k}(t)$ is of the form : 
\begin{equation}\label{eq:form}
{\mit \Delta}^{(2)}_{k,k}(t)\doteq r_1(t)r_2(t)r_3(t)
\end{equation}
where $r_i(t)\doteq r_i(t^{-1})$, $|r_i(1)|=1$ and 
$r_i(t)\doteq_2t^2+t+1$\ ($i=1, 2, 3$).
That is, ${\mit \Delta}^{(2)}_{k,k}(t)$ is decomposed into 
at least three non-trivial factors in $\mathbb{Z}[t, t^{-1}]$.
We set $d_i$ as the degree (span) of $r_i(t)$\ $(i=1, 2, 3)$, and
assume $d_1\le d_2\le d_3$.
There are two cases : 

\medskip

\noindent
{\bf Case 1}\ $k\equiv 0\ (\mathrm{mod}\ \! 4)$.

\medskip

By Lemma \ref{lm:StAlex}, we have the mod $8$ Alexander polynomial : 
$${\mit \Delta}^{(2)}_{k,k}(t)\doteq_8
t^6-t^5+2t^4+3t^3+2t^2-t+1.$$
Since $t^2\pm t+1$ and $t^2\pm 3t+1$ are not mod $8$ divisors of ${\mit \Delta}^{(2)}_{k,k}(t)$,
the case does not occur.

\bigskip

\noindent
{\bf Case 2}\ $k\equiv 2\ (\mathrm{mod}\ \! 4)$.

\medskip

By Lemma \ref{lm:StAlex}, we have the mod $8$ Alexander polynomial : 
$$\begin{array}{rcl}
{\mit \Delta}^{(2)}_{k,k}(t) & \doteq_8 & 
4t^{12}+4t^{11}+3t^9+t^8-2t^7-3t^6-2t^5+t^4+3t^3+4t+4
\medskip\\
& \doteq_8 & 
(t^2-t+1)(4t^{10}+4t^8-t^7+4t^6+3t^5+4t^4-t^3+4t^2+4).
\end{array}$$
We set $s=4t^{10}+4t^8-t^7+4t^6+3t^5+4t^4-t^3+4t^2+4$.
In this case, the $\mathbb{Z}$-degree of ${\mit \Delta}^{(2)}_{k,k}(t)$ is $12$
which is equal to the mod $8$ degree of it.
By the assumption, there are three cases for the triple $(d_1, d_2, d_3)$ : 
$(d_1, d_2, d_3)=(2, 2, 8)$, $(2, 4, 6)$ or $(4, 4, 4)$.
The possibilities of the degree $2$ mod $8$ factors are $t^2\pm t+1$ and $t^2\pm 3t+1$.
Since $t^2\pm t+1$ and $t^2\pm 3t+1$ are not mod $8$ divisors of $s$,
$s$ is decomposed into $s=s_1s_2$ such that
the degrees of $s_1$ and $s_2$ are $4$ and $6$ respectively, 
they are both irreducible, and
$s_1\doteq_2s_2\doteq_2 t^2+t+1$.
By (\ref{eq:form}), $s_1$ and $s_2$ are of the form : 
$$\begin{array}{ccl}
s_1 & \doteq_8 & 2t^4+a_1t^3+a_2t^2+a_1t+2\doteq_2 t^2+t+1,
\medskip\\
s_2 & \doteq_8 & 2t^6+b_1t^5+b_2t^4+b_3t^3+b_2t^2+b_1t+2
\doteq_2 t^2+t+1
\end{array}$$
where $a_1$, $a_2$, $b_2$ and $b_3$ are odd, and $b_1$ is even.
Then the $9$-th coefficient of $s_1s_2$ is odd (non-zero).
However it contradicts the form of $s$.
\qed

\bigskip

At the end of the paper, we raise refined questions realted with Question \ref{qu:cp}:
\begin{qu}\label{qu:mincr}
\begin{enumerate}
\item[(1)]
Is there a prime component-preservingly amphicheiral link
with odd minimal crossing number less than $21$ ?

\item[(2)]
Is there a prime component-preservingly 
$(\varepsilon)$-amphicheiral link
with odd minimal crossing number ?

\end{enumerate}
\end{qu}
About (1), we have already known that there are no such examples
for the case that the minimal crossing number $\le 11$
(cf.\ \cite{Kd3}).
If we need to use an amphicheiral knot 
with odd minimal crossing number, then
the minimal crossing number should be
greater than or equal to $19$ from primeness.
Under the restriction, if there exists
an example $L$ for Question \ref{qu:mincr} (1)
with minimal crossing number $19$, then
$L$ is a $2$-component link such that
(i)\ its components are a knot
with minimal crossing number $15$
and the unknot, 
(ii)\ $\mathrm{lk}\ \! (L)=0$, and
(iii)\ on its diagram realizing the minimal crossing number,
its components are also realizing the minimal crossing numbers
(i.e.\ $15$ and $0$).

\medskip

About (2), our example $L_n$ was
a prime component-preservingly 
$(-, +)$-amphicheiral link with odd minimal crossing number.
In general, the linking number of a $2$-component 
$(\varepsilon)$-amphicheiral link is $0$.
$11_{n247}^2$ in Figure 2 is 
a prime $(\varepsilon)$-amphicheiral link
with odd minimal crossing number.
However it is not component-preservingly 
$(\varepsilon)$-amphicheiral.

\bigskip

{\noindent {\bf Acknowledgements}}\ 
The authors would like to express gratitudes to
Professor Akio Kawauchi,
Professor Taizo Kanenobu, 
Kenji Shibata, members of Topology Seminar at Osaka City University,
and the anonymous referee
for giving them useful comments.

{\footnotesize
\par
\medskip
Teruhisa KADOKAMI\par 
Department of Mathematics,\par
East China Normal University,\par
Dongchuan-lu 500, Shanghai, 200241, China \par
{\tt mshj@math.ecnu.edu.cn}\par
{\tt kadokami2007@yahoo.co.jp}\par
\medskip
Yoji KOBATAKE \par
{\tt koba0726402@gmail.com}\par
}
\end{document}